\documentclass[11pt]{article}
\usepackage[margin=1in]{geometry}
\usepackage[T1]{fontenc}
\usepackage{lmodern}
\usepackage{amsmath,amssymb,amsthm,mathtools}
\usepackage{comfortaa,alphabeta}
\newcommand{\srg}[1]{\text{\normalfont\comfortaa #1}}
\usepackage{microtype,enumitem,booktabs}
\usepackage[dvipsnames]{xcolor}
\usepackage[hidelinks]{hyperref}
\hypersetup{pdftitle={Explicit restricted isometries with exponent 2/3},pdfauthor={Christian Kuemmerle and Dustin G. Mixon}}
\newcommand{\Fp}{\mathbb F_p}

\newcommand{\Fppp}{\mathbb F_{p^3}}
\newcommand{\C}{\mathbb C}

\newcommand{\ep}{e_p}
\newcommand{\one}{\mathbf 1}

\DeclareMathOperator{\rank}{rank}

\newcommand{\ip}[2]{\langle #1,#2\rangle}
\newcommand{\norm}[1]{\lVert #1\rVert}

\newcommand{\DR}{P}
\newtheorem{theorem}{Theorem}
\newtheorem{lemma}[theorem]{Lemma}
\newtheorem{proposition}[theorem]{Proposition}

\theoremstyle{definition}

\newtheorem{problem}[theorem]{Problem}
\theoremstyle{remark}

\title{Explicit restricted isometries with $k=m^{2/3-o(1)}$}
\author{Christian K\"ummerle\thanks{Dr.\ Yvette Kanouff School of Data, Mathematical and Statistical Sciences \& Department of Computer Science, University of Central Florida, Orlando, FL}
\and Dustin G.\ Mixon\thanks{Department of Mathematics \& Translational Data Analytics Institute, The Ohio State University, Columbus, OH}}
\date{}
\begin{document}
\maketitle
\begin{abstract}
We construct an explicit family of real $m\times n$ restricted isometries with $m\to\infty$, $m/n\to0$, $\delta_k\to0$, and $k=m^{2/3-o(1)}$.
We combine Gershgorin’s circle theorem with a flat RIP estimate to reduce the problem to choosing a large vertex set in a strongly regular graph with the property that every small subset spans few edges. 
We select this vertex set using an incidence geometric result due to Zahl. 
The proof only uses the entrywise absolute value of the Gram matrix despite the conventional wisdom that such proof techniques exhibit a barrier of $m^{1/2}$. 
We prove that this barrier is actually $m^{2/3}$, and we show that our particular construction exhibits the same barrier.
\end{abstract}

\section{Introduction}\label{sec:introduction}

Cand\`es and Tao~\cite{CT} introduced the notion of \textit{restricted isometry} to quantify how well a measurement matrix preserves the geometry of sparse vectors in the context of compressed sensing.
For an $m\times n$ matrix $\Phi$, its restricted isometry constant $\delta_k$ is the smallest $\delta$ for which
\[
 (1-\delta)\norm z_2^2\le\norm{\Phi z}_2^2\le(1+\delta)\norm z_2^2
 \qquad\text{whenever }|\operatorname{supp}z|\le k.
\]
Equivalently,
$\delta_k=\max_{|S|\le k}\norm{\Phi_S^*\Phi_S-I}_{2\to2}$, where $\Phi_S$ is the submatrix with columns indexed by $S$.

Random matrices have nearly minimal restricted isometry constants. 
The problem of obtaining comparable explicit constructions was posed by Tao in~\cite{Tao}. 
Here, \textit{explicit} means that the entries can be generated deterministically in time polynomial in the matrix dimensions and the requested output precision. 
We seek an explicit family of $m\times n$ matrices such that
\begin{equation}\label{eq:goal}
m\longrightarrow\infty,
\qquad 
\frac mn\longrightarrow0,
\qquad
\delta_k\longrightarrow0,
\qquad 
k=m^{\alpha-o(1)}
\end{equation}
for the largest possible $\alpha$. 
Random constructions as in~\cite{CTnear,MPT,RV,BDDW} reach an exponent of $\alpha=1$, while elementary explicit constructions only have $\alpha=1/2$.

\begin{table}[t]
\centering
\begin{tabular}{@{\hspace{1em}}l@{\hspace{3em}}l@{\hspace{1em}}}
\toprule
$\alpha$ & Location\\
\midrule
$1/2$ & Thm.~3.1 in \cite{DeVore}\\
$1/2+5.51\cdot10^{-28}$ & Thm.~1 in \cite{BDFKK}\\
$1/2+9.08\cdot10^{-16}$ & Thm.~1 in \cite{BDFKKSTOC}\\
$1/2+8.15\cdot10^{-8}$ & Thm.~1 in \cite{FKS}\\
$2/3$ & Thm.~\ref{thm:main}\\
\bottomrule
\end{tabular}
\caption{Unconditional exponents for known restricted isometries.}
\label{tab:exponents}
\end{table}

To see the source of the exponent $1/2$, suppose the columns $\varphi_1,\ldots,\varphi_n$ of $\Phi$ have unit norm and coherence
\[
\mu:=\max_{i\ne j}|\ip{\varphi_i}{\varphi_j}|.
\]
Gershgorin's circle theorem gives
\begin{equation}\label{eq:coherence}
 \delta_k\le(k-1)\mu,
\end{equation}
suggesting we design $\Phi$ to have incoherent columns. 
This can be accomplished, for example, using polynomial phase vectors over finite fields.
For any odd prime $p$, define $e_p\colon\Fp\to\C$ by
\[
\ep(t)
=\exp(2\pi i t/p).
\]
Then the unit vectors $v_f:=p^{-1/2}(\ep(f(x)))_{x\in\Fp}$ are incoherent: if $r:=\deg(f-g)$ is positive and smaller than $p$, then Weil's bound \cite{Weil} gives
\[
 |\ip{v_f}{v_g}|\le\frac{r-1}{\sqrt p}.
\]
If we take all bounded-degree polynomials $f$ satisfying $f(0)=0$, then the vectors $v_f$ have coherence $O(p^{-1/2})$ for all sufficiently large $p$. 
Assuming the degree bound is at least $2$, then $m=p$ satisfies $m/n\to0$ as in~\eqref{eq:goal}, while \eqref{eq:coherence} implies $\delta_k\to0$ for $k=m^{1/2}(\log m)^{-1}$ (say), giving $\alpha=1/2$.
For a fundamental limit to this approach, the Welch bound~\cite{Welch} states that
\begin{equation}\label{eq.welch}
 \mu\ge\sqrt{\frac{n-m}{m(n-1)}}
\end{equation}
for any collection of $n$ unit vectors in $\C^m$, and this lower bound is asymptotically $m^{-1/2}$ when $m/n\to0$.
Thus, if also $m\to \infty$, the  estimate~\eqref{eq:coherence} can certify $\delta_k\to0$ only when $k=o(m^{1/2})$.

A different approach is necessary to certify an exponent beyond $1/2$.
Since the Gershgorin estimate~\eqref{eq:coherence} ignores the phases in the Gram matrix $\Phi^*\Phi$, one possible improvement is to account for \textit{phase cancelation}.
This was first accomplished by Bourgain, Dilworth, Ford, Konyagin, and Kutzarova in~\cite{BDFKK}, who used a particular family of quadratic polynomials:
\[
 u_{a,b}(x)=p^{-1/2}\ep(ax^2+bx),
 \qquad (a,b)\in\mathcal A\times\mathcal B.
\]
For distinct column indices, if $a=a'$, then $u_{a,b}$ is orthogonal to $u_{a',b'}$, while if $a\neq a'$, then the modulus of $\ip{u_{a,b}}{u_{a',b'}}$ is $p^{-1/2}$, and its phase is determined by a Gauss sum.
Bourgain et al.\ convert the desired cancelations into questions about the additive combinatorics of $\mathcal A$ and $\mathcal B$, and then a careful choice of these sets delivers an exponent strictly larger than $1/2$. 
See~\cite{Mixon} for a friendly presentation of these ideas and~\cite{BDFKKSTOC,FKS} for subsequent refinements.
Table~\ref{tab:exponents} summarizes these results.

To obtain restricted isometry from phase cancelation, one establishes \emph{flat RIP} by bounding
\begin{equation}\label{eq:intro-flat}
 \bigg|\bigg\langle\sum_{i\in I}\varphi_i,\sum_{j\in J}\varphi_j\bigg\rangle\bigg|
 \le\theta\sqrt{|I||J|}
\end{equation}
for disjoint sets $I,J$ of size at most $k$. 
A multiscale argument then gives $\delta_k\lesssim\theta\log(2k)$ for unit-norm columns; see~\cite{BDFKK,BFMW}. 
These normalized rectangular sums are the disjoint-set version of the \emph{discrete Rayleigh norm} studied by Lev~\cite{Lev}. 
We make this connection precise in Section~\ref{sec:motivation}.

One can also leverage flat RIP to obtain conditional constructions of explicit restricted isometries. 
Bandeira, Mixon, and Moreira~\cite{BMM} show that a conjectured discrepancy bound for the Paley graph implies restricted isometry beyond the square-root threshold for the Paley equiangular tight frame. 
Satake~\cite{Satake} extends this conditional approach to arbitrary odd primes. 
These Paley constructions have row-to-column ratio $1/2$ rather than the vanishing ratio required in~\eqref{eq:goal}, but they illustrate how number-theoretic pseudorandomness can supply flat RIP estimates.

In this paper, we use a different approach: we apply the induced $1$-norm underlying Gershgorin to the small entries of the Gram matrix $\Phi^*\Phi$, and the discrete Rayleigh norm underlying flat RIP to its large entries. 
For the latter estimate, we take entrywise absolute values and bound the density of every small induced subgraph of the support graph. 
In particular, \textit{our proof does not use phase cancelation}, and yet we obtain $\alpha=2/3$ unconditionally.

\begin{theorem}
\label{thm:main}
There is an explicit family of real $m\times n$ matrices with unit-norm columns such that
\[
m\longrightarrow\infty,
\qquad 
\frac mn\longrightarrow0,
\qquad
\delta_k\longrightarrow0,
\qquad 
k=m^{2/3-o(1)}.
\]
\end{theorem}

Section~\ref{sec:motivation} develops the hybrid norm estimate, reduces it to the density of subgraphs, and proves a $2/3$~barrier for arguments that do not use phase cancelation. 
Section~\ref{sec:frames} uses this barrier to motivate the use of two-distance tight frames and identifies a strongly regular graph whose dense subgraphs were studied by Zahl in~\cite{Zahl}. 
Section~\ref{sec:proof} then uses this information to carefully select vertices from this graph and prove Theorem~\ref{thm:main}. 
Finally, Section~\ref{sec:barrier} proves that our restricted isometry fails to pass the $2/3$~barrier, even when accounting for phase cancelation.

\section{A new way to estimate restricted isometry constants}\label{sec:motivation}

In this section, we combine well-known estimates in a new way that we later use to break the square-root bottleneck $k\asymp m^{1/2}$ while ignoring phase cancelation in the Gram matrix.

\subsection{Expressing Gershgorin and flat RIP in terms of matrix norms}

The induced $1$-norm of a matrix $M$ is given by
\[
\norm M_{1\to1}
:=\max_j\sum_i|M_{ij}|.
\]
If $M$ is self-adjoint, the spectral norm $\norm M_{2\to2}$ can be bounded by Gershgorin's circle theorem, placing every eigenvalue in a disc $\{z:|z-M_{jj}|\le\sum_{i\ne j}|M_{ij}|\}$, which implies
\begin{equation*}
\norm M_{2\to2}\le\norm M_{1\to1},
\end{equation*}
in matrix-norm form.
To similarly capture flat RIP, we  consider the \emph{discrete Rayleigh norm}~\cite{Lev} defined by
\[
\norm M_{\DR}
:=
 \max_{\varnothing\ne I,J\subseteq[k]}
 \frac{|\one_I^*M\one_J|}{\sqrt{|I||J|}}.
\]
Then Theorem~2 in~\cite{Lev} implies the bound
\begin{equation}\label{eq:Lev}
\norm M_{2\to2}
\lesssim\log(2k)\norm M_{\DR}.
\end{equation}
By the following lemma, this is a matrix-norm form of the flat RIP estimate.

\begin{lemma}
\label{lem:flat-rip-dr}
Let $\Phi\in\mathbb{C}^{m\times n}$ have
unit-norm columns $\varphi_1,\ldots,\varphi_n\in\mathbb{C}^m$. 
Let $\theta_k(\Phi)$ denote the
smallest $\theta \geq 0 $ for which \eqref{eq:intro-flat} holds for all disjoint $I,J\subseteq[n]$ with $|I|,|J|\le k$. 
Define
\[
\rho_k(\Phi):=
  \max_{\substack{\varnothing\ne S\subseteq[n]\\ |S|\le k}}
  \bigl\|\Phi_S^*\Phi_S-I\bigr\|_{\DR},
\]
where $\Phi_S$ consists of the columns of $\Phi$ indexed by $S$. 
Then
\[
  \theta_k(\Phi)\le \rho_{2k}(\Phi)
  \le 2\theta_k(\Phi).
\]
\end{lemma}

\begin{proof}
For $\theta_k(\Phi)\le \rho_{2k}(\Phi)$, the case $\theta_k(\Phi)=0$ is immediate. Otherwise, since there are only finitely many admissible pairs, select disjoint nonempty $I,J\subseteq[n]$ with $|I|,|J|\leq k$ such that
\[
\bigg|\bigg\langle\sum_{i\in I}\varphi_i,\sum_{j\in J}\varphi_j\bigg\rangle\bigg|
=\theta_k(\Phi)\sqrt{|I||J|}.
\]
Then
\[
\theta_k(\Phi)
=\frac{|\one_I^*(\Phi_{I\cup J}^*\Phi_{I\cup J}-I)\one_J|}{\sqrt{|I||J|}}
\leq\bigl\|\Phi_{I\cup J}^*\Phi_{I\cup J}-I\bigr\|_{\DR}
\leq\rho_{2k}(\Phi).
\]
For $\rho_{2k}(\Phi)\leq 2\theta_k(\Phi)$, fix any nonempty $S\subseteq[n]$ with $|S|\le 2k$ and any nonempty $I,J\subseteq S$. Adjoin auxiliary labels to $S$ to form a set $T$ of size $2k$.
Choose a uniformly random $k$-element subset $R$ of $T$, and put $A=S\cap R$ and $B=S\setminus R$, which are disjoint and have at most $k$ elements each.
Then for distinct $i,j\in S$, we have
\[
  \mathbb{P}(i\in A,\, j\in B)=\frac12\cdot\frac{k}{2k-1},
\]
as well as $\mathbb E|I\cap A|=\frac{1}{2}|I|$ and $\mathbb E|J\cap B|=\frac{1}{2}|J|$. 
Since $\Phi_S^*\Phi_S-I$ has zero diagonal by assumption, it follows that
\begin{align*}
  \frac{k}{2(2k-1)}|\one_I^*(\Phi_S^*\Phi_S-I)\one_J|
  &=\bigl|\mathbb{E}\,\one_{I\cap A}^*(\Phi_S^*\Phi_S-I)\one_{J\cap B}\bigr|\\
  &\le \theta_k(\Phi)\,
       \mathbb{E}\bigg[\sqrt{|I\cap A|}\sqrt{|J\cap B|}\bigg]\\
  &\le \theta_k(\Phi)\,
       \sqrt{\mathbb{E}|I\cap A|}\sqrt{\mathbb{E}|J\cap B|}
  =\frac{\theta_k(\Phi)}{2}\sqrt{|I||J|},
\end{align*}
where the first inequality follows from the triangle inequality and the definition of $\theta_k(\Phi)$, and the second follows from Cauchy--Schwarz.
Rearranging and maximizing over $I$, $J$, and $S$ then gives 
\[
\rho_{2k}(\Phi)
\leq\frac{2k-1}{k}\theta_k(\Phi)
\leq 2\theta_k(\Phi).
\qedhere
\]
\end{proof}

\subsection{Defying conventional wisdom}

To date, the community has used the following program to identify explicit restricted isometries:
\begin{center}
Find unit vectors of low coherence with estimable phase cancelation in the Gram matrix.
\end{center}
By \textit{low coherence}, we mean coherence on the order of $m^{-1/2+o(1)}$.
There are many reasons to restrict the search to such ensembles.
First, unit vectors have restricted isometry constant $\delta_1=0$, and then their coherence equals $\delta_2$.
Second, the Gershgorin estimate further ensures that $\delta_k$ is small for $k=m^{1/2-o(1)}$.
Third, the random constructions achieving sparsity~$k=m^{\alpha-o(1)}$ with exponent~$\alpha=1$ also exhibit low coherence (see, e.g.,~\cite{CaiJiang}).
However, in practice, we encounter difficulty by restricting to low-coherence ensembles.
First of all, explicit examples of such ensembles are difficult to come by~\cite{FMtables}.
Second, by minimizing coherence, we frequently encounter linearly dependent subsets of size $O(m^{1/2})$ as in~\cite{FMT,FJKM}.
For some reason, there appears to be a tradeoff between minimizing $\delta_2$ and maximizing $k$ for which $\delta_k$ is small (cf.~\cite{KingSpread}).

This suggests a different approach that defies conventional wisdom:
\begin{center}
Consider unit vectors that \textit{do not} exhibit low coherence.
\end{center}
Let $M$ denote a $k\times k$ principal submatrix of $\Phi^*\Phi-I$.
In the absence of low coherence, the Gershgorin estimate is weak.
Let $Y$ denote the version of $M$ that has entries zeroed out if they are too large for Gershgorin, and put $X:=M-Y$.
Then, applying Lev's bound~\eqref{eq:Lev} to~$X$, we have
\begin{equation}
\label{eq.split}
\|M\|_{2\to2}
\leq \|X\|_{2\to2}+\|Y\|_{2\to2}
\lesssim\log(2k)\|X\|_{\DR}+\|Y\|_{1\to1}.
\end{equation}
By giving the large entries of $M$ to the discrete Rayleigh norm, we avoid a failure mode of the Gershgorin estimate, but we also force ourselves to interact with the flat RIP estimate.
While this estimate allows us to treat phase cancellation, it is extraordinarily technical and cumbersome.
Since we are already defying conventional wisdom, we might as well go all the way:
\begin{center}
\textit{Do not} estimate phase cancelation in the Gram matrix.
\end{center}
We can avoid working with phases in $X$ by taking the entrywise absolute value, since the triangle inequality gives the bound
\begin{equation}
\label{eq.dr abs}
\|X\|_{\DR}
\leq\big\||X|\big\|_{\DR}.
\end{equation}
Going further, if we let $A$ denote the binary matrix that indicates the support of $X$, then 
\begin{equation}
\label{eq.pass to graph}
\big\||X|\big\|_{\DR}
\leq\|X\|_{\max}\|A\|_{\DR}
\leq\mu\|A\|_{\DR},
\end{equation}
where again, $\mu$ is the coherence of $\Phi$.
We may interpret $A$ as the adjacency matrix of a graph, in which case bounding $\|A\|_{\DR}$ reduces to showing that all of its subgraphs have low density; this is made rigorous in the following lemma, which is implicit in the proof of Lemma~2.6 in~\cite{CZZ}.

\begin{lemma}
\label{lem:density-norm}
Given a finite simple graph with vertex set $V$, let $e(U)$ denote the number of edges induced by $U\subseteq V$. 
Suppose that, for some $c\ge 0$ and $1\le\beta\le 2$,
\[
e(U)\le c|U|^\beta
\]
for every $U\subseteq V$.
Then the adjacency matrix $A$ of the graph satisfies
\[
\|A\|_{\DR}
\le 2^{\beta+1}c\,|V|^{\max\{1/2,\,\beta-1\}}.
\]
\end{lemma}

\begin{proof}
Fix nonempty $I,J\subseteq V$.
We may assume $|I|\le |J|$. 
Choose $J'$ uniformly among the $|I|$-element subsets of $J$. 
Each vertex of $J$ belongs to $J'$
with probability $|I|/|J|$, and so
\[
    \frac{|I|}{|J|}\one_I^*A\one_J
    = \mathbb{E}\,\one_I^*A\one_{J'}
    \le 2\,\mathbb{E}\,e(I\cup J')
    \le 2c(2|I|)^\beta.
\]
The first inequality applies the fact that $\one_I^*A\one_{J'}$ counts each edge induced by $I\cup J'$ at most twice, while the second inequality uses $|I\cup J'|\le 2|I|$.
Rearranging gives
\[
    \frac{\one_I^*A\one_J}{\sqrt{|I||J|}}
    \le 2^{\beta+1}c\,|J|^{1/2}\,|I|^{\beta-3/2}.
\]
If $\beta\le 3/2$, then $|I|^{\beta-3/2}\le 1$, giving the bound $2^{\beta+1}c|V|^{1/2}$. 
If $\beta\ge 3/2$, then $|I|,|J|\le |V|$ gives the bound $2^{\beta+1}c|V|^{\beta-1}$. Maximizing over
$I$ and $J$ proves the claim.
\end{proof}

Overall, suppose $\varphi_1,\ldots,\varphi_n\in\mathbb{C}^m$ are unit norm with coherence $\mu$, select a threshold $\nu\in(0,\mu)$, and define a graph in which $i,j\in[n]$ with $i\neq j$ are adjacent when $|\langle\varphi_i,\varphi_j\rangle|>\nu$.
Furthermore, suppose every subset $U\subseteq[n]$ of size at most $k$ induces only $e(U)\lesssim |U|^\beta$ edges in this graph.
Let $A$ denote the adjacency matrix of this graph, and let $H$ denote the version of $\Phi^*\Phi-I$ with entries of magnitude greater than $\nu$ zeroed out.
We combine the estimates~\eqref{eq.split}, \eqref{eq.dr abs}, and~\eqref{eq.pass to graph}, and then apply Lemma~\ref{lem:density-norm} to obtain
\begin{align}
\nonumber
\delta_k
&=\max_{|S|\leq k}\|\Phi_S^*\Phi_S-I\|_{2\to2}\\
\nonumber
&\lesssim \mu\log(2k)\max_{|S|\leq k}\| A_S\|_{\DR}+\max_{|S|\leq k}\|H_S\|_{1\to1}\\
\label{eq.idea of bound}
&\lesssim k^{\max\{1/2,\beta-1\}}\log(2k)\cdot\mu +k\cdot\nu.
\end{align}
To allow for the largest possible $\mu$, we take the smallest possible $\beta$, though there is no benefit in taking $\beta$ smaller than $3/2$, so we target $\beta=3/2$.
This allows us to accommodate coherence $\mu$ as large as $k^{-1/2-o(1)}$, which dwarfs the Welch bound~\eqref{eq.welch} when $k\ll m$.
Meanwhile, the threshold $\nu$ will have to be $o(k^{-1})$, which is much smaller than the Welch bound in the desired regime $k\gg m^{1/2}$.

\subsection{The limit of phase agnosticism}\label{subsec:absolute-barrier}

Notice that taking the entrywise absolute value of a matrix does not decrease its spectral norm.
Indeed, if the unit vector $x$ is a top right singular vector of $M$, then
\[
\|M\|_{2\to2}^2
=\|Mx\|^2
=\sum_i \bigg|\sum_j M_{ij}x_j\bigg|^2
\leq\sum_i \bigg|\sum_j |M_{ij}||x_j|\bigg|^2
=\big\||M||x|\big\|^2
\leq \big\||M|\big\|_{2\to2}^2.
\]
It follows that
\[
\Delta_k
:=\max_{|S|\le k}\big\||\Phi_S^*\Phi_S-I|\big\|_{2\to2}
\]
is an upper bound on the restricted isometry constant $\delta_k$.
Notice that both the Gershgorin estimate and the estimate obtained by combining \eqref{eq.split} and \eqref{eq.dr abs} in the previous subsection are actually upper bounds on $\Delta_k$.
The following result reveals the limit of any such technique.

\begin{theorem}
\label{thm:absolute-barrier}
Suppose $\Phi\in\C^{m\times n}$ has unit-norm columns and $2\leq k\leq m\leq n$.
Then
\begin{equation}\label{eq:barrier}
 \Delta_k^2\ge
 \frac{(n-m)(k-1)^{3/2}}
 {m\bigl(n-1+(k-1)^{3/2}\bigr)}.
\end{equation}
In particular, if $m \to \infty$ and $m/n\to0$, then $\Delta_k\to0$ requires $k=o(m^{2/3})$.
\end{theorem}

\begin{proof}
Put $H:=|\Phi^*\Phi-I|$. 
We first apply Welch's bound on the frame potential, namely, that for any $\Phi\in\mathbb{C}^{m\times n}$ with unit-norm columns,  $\|\Phi^*\Phi\|_F^2\geq n^2/m$.
Relating this to $H$ gives
\begin{equation}
\label{eq.frame potential}
\|H\|_F^2
=\|\Phi^*\Phi\|_F^2-n
\geq\frac{n^2}{m}-n.
\end{equation}
Next, we find an upper bound on $\|H\|_F^2$ in terms of $\Delta_k$.
In each row of $H$, order the off-diagonal entries as $h_{i,1}\ge\cdots\ge h_{i,n-1}$.
Let $S$ consist of $i$ and the $k-1$ column indices corresponding to $h_{i,1},\ldots,h_{i,k-1}$, and let $H_S$ denote the principal submatrix of $H$ indexed by $S$.
Then
\[
\sum_{j=1}^{k-1}h_{i,j}^2
=\|H_Se_i\|^2
\leq\|H_S\|_{2\to2}^2
\le\Delta_k^2.
\]
For every $j\geq k$, it holds that
\[
h_{i,j}^2
\leq h_{i,k-1}^2
\leq \frac{1}{k-1}\sum_{j=1}^{k-1}h_{i,j}^2
\leq \frac{\Delta_k^2}{k-1}.
\]
Taking the square root of both sides above gives the second inequality below:
\[
\sum_{j=1}^{n-1}h_{i,j}^2
\leq \sum_{j=1}^{k-1}h_{i,j}^2+\sum_{j=k}^{n-1}h_{i,j}^2
\leq \Delta_k^2+\frac{\Delta_k}{\sqrt{k-1}}\sum_{j=1}^{n-1}h_{i,j}.
\]
Summing over $i$ then gives
\begin{equation}
\label{eq.bound by entrywise 1 norm}
\|H\|_F^2
\leq n\Delta_k^2+\frac{\Delta_k}{\sqrt{k-1}}\|H\|_{1,1},
\end{equation}
where $\|H\|_{1,1}:=\sum_{i,j=1}^n|H_{ij}|$ denotes the entrywise $\ell_1$-norm,  which equals $\sum_{i,j=1}^n H_{ij}$ due to its non-negative entries. Now we seek an upper bound on $\|H\|_{1,1}$.
For every $S\subseteq[n]$ of size $k$, we have
\[
\frac{1}{k}\sum_{i,j\in S}H_{ij}
=\frac{\one^*H_S\one}{k}
\leq\|H_S\|_{2\to2}
\leq\Delta_k.
\]
Next, we take the average, observing that $H_{ii}=0$ and that each $\{i,j\}$ with $i\neq j$ is contained in $\binom{n-2}{k-2}$ different sets of size $k$:
\[
\Delta_k
\geq\frac{1}{\binom{n}{k}}\sum_{\substack{S\subseteq[n]\\|S|=k}}\frac{1}{k}\sum_{i,j\in S}H_{ij}
=\frac{\binom{n-2}{k-2}}{k\binom{n}{k}}\sum_{i,j\in[n]}H_{ij}
=\frac{k-1}{n(n-1)}\|H\|_{1,1}.
\]
Combining with~\eqref{eq.frame potential} and~\eqref{eq.bound by entrywise 1 norm} then gives
\[
\frac{n^2}{m}-n
\leq n\Delta_k^2+\frac{\Delta_k}{\sqrt{k-1}}\cdot\frac{n(n-1)}{k-1}\Delta_k,
\]
and rearranging gives~\eqref{eq:barrier}.
For the final claim, suppose $m/n\to0$ and $\Delta_k\to0$.
Then~\eqref{eq:barrier} implies
\[
\frac{(k-1)^{3/2}}{n-1+(k-1)^{3/2}}
\leq\frac{m}{n-m}\Delta_k^2\to0,
\]
and so $(k - 1)^{3/2}=o(n)$.
It follows that
\[
\Delta_k^2
\geq\frac{(n-m)(k-1)^{3/2}}{m\bigl(n-1+(k-1)^{3/2}\bigr)}
=\frac{n-m}{n-1+(k-1)^{3/2}}\cdot\frac{(k-1)^{3/2}}{m}
=(1+ o(1)) \frac{(k-1)^{3/2}}{m}.
\]
Thus, $\Delta_k\to 0$ forces $(k-1)^{3/2}/m \to 0$, i.e., $k = o(m^{2/3})$.
\end{proof}

Previously, the only phase-agnostic approach factored through the Gershgorin estimate.
Consider the \textit{cumulative coherence function}, introduced by Tropp in~\cite{TroppGreed}:
\[
\mu_1(k-1)
:=\max_{|S|\leq k}\|\Phi_S^*\Phi_S-I\|_{1\to1}.
\]
In~\cite{SV}, Schnass and Vandergheynst proved that
\[
\mu_1(k-1)
\ge (k-1)\sqrt{\frac{n-m}{m(n-1)}}
\]
whenever $(k-1)^2<n-1$.
While it holds that $\mu_1(k-1)\geq\Delta_k$, in the setting where $m\to\infty$ and $m/n\to0$, the vanishing $\mu_1(k-1)\to0$ requires $k=o(m^{1/2})$.
This phenomenon has become known as the \textit{square-root bottleneck}~\cite{TroppConditioning}.
To quote Rauhut from~\cite{Rauhut}, the Gershgorin estimate ``seems to be the optimal estimate one may obtain by taking into account only the absolute values of the Gramian matrix.''
As we will see, it is possible to approach equality in Theorem~\ref{thm:absolute-barrier}, and so this square-root bottleneck is actually an artifact of the Gershgorin estimate rather than of phase agnosticism.

\section{Two-distance tight frames and strongly regular graphs}\label{sec:frames}

Theorem~\ref{thm:absolute-barrier} suggests aiming for sparsity $k=m^{2/3-o(1)}$.
At this scale, the two terms in~\eqref{eq.idea of bound} with $\beta=3/2$ suggest looking for large inner products of magnitude $\mu\asymp m^{-1/3}$ and small inner products of magnitude $\nu\asymp m^{-2/3}$. 
In order to approach the bound in Theorem~\ref{thm:absolute-barrier}, we also want equality in the frame potential bound~\eqref{eq.frame potential}. 
This motivates the use of \emph{two-distance tight frames}~\cite{BGOY}, which are matrices $\Phi\in\mathbb{R}^{m\times n}$ for which there exist $a,b,c\in\mathbb{R}$ with $a\neq b$ such that
\begin{itemize}
\item[(a)]
the diagonal entries of $\Phi^*\Phi$ all equal $1$,
\item[(b)]
the off-diagonal entries of $\Phi^*\Phi$ all reside in $\{a,b\}$, and
\item[(c)]
$\Phi$ forms a $c$-tight frame, meaning $\Phi\Phi^*=cI$.
\end{itemize}
Note that (c) is equivalent to $c^{-1}\Phi^*\Phi$ being an orthogonal projection matrix of rank $m$, and so $c$ necessarily equals $n/m$.
In our use case, $a$ and $b$ will have absolute value $\mu$ and $\nu$, respectively.
We will construct such a frame with many more columns than rows so that we can select a subframe whose Gram matrix determines the low-density subgraphs necessary for our bound~\eqref{eq.idea of bound}.

A graph is said to be \textit{strongly regular} with parameters $(\srg{v},\srg{k},\srg{\lambda},\srg{\mu})$ if it has exactly $\srg{v}$~vertices, each of degree~$\srg{k}$, such that adjacent vertices have $\srg{\lambda}$ common neighbors and distinct nonadjacent vertices have $\srg{\mu}$ common neighbors; see~\cite{BrouwerVanMaldeghem} for a standard reference. 
(We use a sans-serif font for these parameters to thwart a notational collision with sparsity~$k$ and coherence~$\mu$.)
To avoid degeneracies, it is standard to assume that the graph and its complement are both connected.
Write $A$ for the adjacency matrix, $J =\one \one^{\top}$ for the all-ones matrix, and $\bar A=J-I-A$ for the adjacency matrix of the complementary graph.

Strongly regular graphs have a rich spectral theory.
Since the graph is $\srg{k}$-regular, we have $A\one=\srg{k}\one$.
Next, counting walks of length two gives
\[
A^2
=\srg{k}I+\srg{\lambda}A+\srg{\mu}\bar{A}
=(\srg{k}-\srg{\mu})I+(\srg{\lambda}-\srg{\mu})A+\srg{\mu}J.
\]
If we restrict the above equation to the orthogonal complement of the all-ones vector, then $A$ satisfies a polynomial of degree~$2$ on this subspace. 
It follows that the adjacency matrix has two eigenvalues on this subspace, and furthermore, one may show that the orthogonal projection $P$ onto either eigenspace belongs to the span of $I,A,J$. 
Writing $m=\rank P$, we see that $P$ has constant diagonal~$m/\srg{v}$.
It follows that $G:=(\srg{v}/m)P$ is the Gram matrix of a two-distance tight frame:
\[
G=I\pm\mu A\mp \nu\bar{A}.
\]
Indeed, the signs on $A$ and $\bar{A}$ cannot both be positive since $G$ is a scaled projection orthogonal to~$\one$, and they cannot both be negative since then $G$ would be the Laplacian of a weighted connected graph, meaning the kernel is too small (i.e., just the span of $\one$).
To obtain a restricted isometry with many more columns than rows, we seek an eigenspace with $m=o(\srg{v})$.
It turns out that our requirements force some relatively strong scaling laws between the strongly regular graph parameters:

\begin{lemma}\label{lem:srg-scaling}
Consider a sequence of strongly regular graphs with varying parameters
$(\srg{v},\srg{k},\srg{\lambda},\srg{\mu})$. 
For each of these graphs, let $P$ denote the orthogonal projection onto an $m$-dimensional eigenspace of its adjacency matrix $A$. 
Suppose
\[
    G:=\frac{\srg{v}}{m}P=I+\varepsilon(\mu A-\nu\bar A),
\]
for some $\varepsilon\in\{\pm1\}$ and $\mu>\nu>0$. If $m\to\infty$, $m=o(\srg{v})$, and
$\mu\asymp m^{-1/3}$, then
\[
\nu\asymp m^{-2/3},
\qquad
\frac{\srg{k}}{\srg{v}}
\sim\frac{\nu}{\mu}
\sim\frac{1}{m\mu^2}
\asymp m^{-1/3},
\qquad
\srg{\mu}
\sim\frac{\srg{k}^2}{\srg{v}},
\qquad
\frac{\srg{\lambda}}{\srg{\mu}}
=1+\varepsilon m\mu^3+o(1).
\]
In particular, $\srg{\lambda}=o(\srg{\mu})$ if and only if
$\varepsilon=-1$ eventually and $\mu\sim m^{-1/3}$, in which case $\nu\sim m^{-2/3}$.
\end{lemma}

See the appendix for a proof.
These scaling laws restrict our attention to a handful of known infinite families of strongly regular graphs.
To select from this handful, recall that we intend to pass to a subgraph with the property that all sufficiently small subgraphs have few edges.
We will accomplish this by applying a theorem of Zahl~\cite{Zahl}, which analyzes the high-density subgraphs of the \textit{elliptic affine polar graph} $\mathrm{VO}_4^-(p)$ for a prime $p\equiv 3\bmod 4$.
This graph has vertex set $\mathbb F_p^4$, and two vertices are adjacent when their difference is isotropic under some fixed nondegenerate quadratic form of Witt index one~\cite{Orel}; see the next section for an explicit description.
The strong regularity of this graph is classical~\cite[p.~375, family C.12]{Hubaut}.
Its parameters are
\[
(\srg{v},\srg{k},\srg{\lambda},\srg{\mu})
=\bigl(p^4,(p-1)(p^2+1),p-2,p(p-1)\bigr)
\sim(p^4,p^3,p,p^2),
\]
which notably satisfy $\srg{\lambda}=o(\srg{\mu})$.
For this graph, the smaller nontrivial eigenspace gives
\begin{equation}\label{eq.params}
m=(p-1)(p^2+1)\sim p^3,
\qquad
\mu=\frac{p^2-p+1}{(p-1)(p^2+1)}\sim m^{-1/3},
\qquad
\nu=\frac{1}{p^2+1}\sim m^{-2/3}.
\end{equation}
In particular, these parameters satisfy the scaling laws in Lemma~\ref{lem:srg-scaling}.
In the next section, we apply the main theorem from~\cite{Zahl} to identify a large subgraph of this strongly regular graph that delivers an explicit restricted isometry.

\section{Selecting a subframe and proving restricted isometry}\label{sec:proof}

Fix a prime $p\equiv3\bmod4$, and consider the $(p-1)(p^2+1)\times p^4$ two-distance tight frame derived from the strongly regular graph $\mathrm{VO}_4^-(p)$, as described in Section~\ref{sec:frames}.
We will select $n=Lp^3$ columns of this frame to form a restricted isometry.
Since the aspect ratio is $m/n\sim1/L$, we may let $L$ tend to infinity as slowly as we please. 
To apply~\eqref{eq.idea of bound} with $\beta=3/2$, we seek a set $V'$ of $n$ vertices from the graph $\mathrm{VO}_4^-(p)$ with the property that
\begin{equation}
\label{eq:edge-goal}
    e(U)\lesssim L|U|^{3/2}
\end{equation}
for every $U\subseteq V'$ of size at most $p^2$.

\subsection{The configurations to avoid}

Recall that our vertex set is $\mathbb{F}_p^4$ for some prime $p\equiv 3\bmod 4$, and we say distinct points $a,b\in\mathbb{F}_p^4$ are adjacent if $Q(a-b)=0$, where $Q$ is a nondegenerate quadratic form of Witt index one.
For an explicit choice of form, we may take\footnote{In the context of Zahl's theorem~\cite{Zahl}, we take $p\equiv 3\bmod 4$ to ensure that $-1$ is not a square in~$\mathbb F_p$.
Also, the vertex $(x,y,z,r)$ determines an oriented sphere with center $(x,y,z)$ and signed radius~$r$, and the adjacency $Q(a-b)=0$ is precisely the contact relation found in equation~(1.1) of~\cite{Zahl}.}
\[
Q(x,y,z,r)
:=x^2+y^2+z^2-r^2.
\]
In order to obtain the estimate~\eqref{eq:edge-goal}, we need to avoid high-density subgraphs.
For example, sizable cliques would be problematic.
In $\mathrm{VO}_4^-(p)$, the maximal cliques are the affine lines with isotropic direction.
Such cliques are known as \textit{pencils}.
We should also avoid any sizable complete bipartite subgraph.
In $\mathrm{VO}_4^-(p)$, the maximal complete bipartite subgraphs with independent sides, each with at least three vertices, are the \emph{complementary conics} of~\cite[Definition~1.6]{Zahl}. 
By~\cite[Remark~1.7]{Zahl}, each side of the bipartite subgraph is an affine plane conic, meaning it is defined by one quadratic equation in an affine two-dimensional plane.
For example, the vertex sets
\[
A:=\{(x,y,0,0):x^2+y^2=1\},
\qquad
B:=\{(0,0,z,r):z^2-r^2=-1\}
\]
are complementary conics, and one may verify that $Q(a-b)=0$ for every $a\in A$ and $b\in B$, and that no other pairs of vertices in $A\cup B$ are adjacent.

Perhaps surprisingly, these are the only configurations we need to avoid.
To see this, let $\mathcal P(U)$ denote the set of pencils containing at least two vertices of $U$. 
Since every edge lies on a unique pencil, it holds that
\[
e(U)
=\sum_{P\in\mathcal P(U)}\binom{|U\cap P|}{2}.
\]
We can control each summand by limiting the number of vertices we
select on each pencil.
To control the number of summands, the following theorem gives that for
$|U|\le p^2$, complementary conics are the only obstruction to a bound of the form $|\mathcal P(U)|\lesssim |U|^{3/2}$.

\begin{theorem}[special case of~{\cite[Theorem~1.8]{Zahl}}\footnote{Out of an abundance of caution, we formalized the proof of Theorem~\ref{thm:zahl} in Lean, along with the theorem of Guth and Zahl~\cite[Theorem~3.8]{GuthZahl} on which it rests.
As is typical in formalization, this required a number of local corrections to the published arguments, but none of them affected the statement of Theorem~\ref{thm:zahl}.}]
\label{thm:zahl}
There is an absolute constant $C$ for which the following holds. Let $p\equiv3\bmod4$ be prime. For each $U\subseteq \mathbb{F}_p^4$ of size $s\le p^2$, at least one of the following holds:
\begin{itemize}
\item[(a)]
The number of pencils containing at least two points from $U$ is at most $C s^{3/2}$.
\item[(b)]
There is a pair of complementary conics, each containing at least $\sqrt{s}$ points from $U$.
\end{itemize}
\end{theorem}

Overall, to obtain~\eqref{eq:edge-goal}, it suffices to bound the number
of points from $V'$ that reside in every pencil and every affine plane conic.

\subsection{Selecting vertices}

Fix $T\subseteq\Fp$ of size $L\ge1$.
To select $Lp^3$ vertices, we will choose $L$ different completions of each $(x,y,z)\in\Fp^3$ to a vertex $(x,y,z,r)\in\Fp^4$. 
If we take each $r$ to be $f(x,y,z)+t$ for $t\in T$, then the selected vertices form $L$ disjoint translates of the graph $r=f(x,y,z)$. 
Since the set of affine lines and affine plane conics is closed under translation, it suffices to select $f$ so that its graph meets every such set in boundedly many points.

Since $\Fp^3$ is finite, $f$ can be implemented as a polynomial.
How complicated does $f$ need to be?
If $f$'s restriction to an affine line in $\Fp^3$ is a polynomial of degree at most $2$, then all $p$ points in the graph's restriction reside on a single line or conic.
Thus, we want $f$ to be at least cubic on every affine line.
To accomplish this, we identify\footnote{Any $\Fp$-linear identification would do the trick, but for concreteness, we identify $(u,v,w)\in\Fp^3$ with $u+vX+wX^2\in\Fp[X]/(f)$, where $f=X^3+aX^2+bX+c$ is the irreducible cubic whose coefficient triple $(a,b,c)\in\{0,\ldots,p-1\}^3$ comes first in lexicographic order.} $\Fp^3$ with $\Fppp$ and take $f$ to be the field norm $N\colon\Fppp\to\Fp$. 
Importantly, the coefficient of $s^3$ in $N(c+sv)$ is $N(v)$, which is nonzero
whenever $v$ is nonzero. 

Overall, we take
\begin{equation}
\label{eq:selection}
V'_T
:=\{(c,N(c)+t):c\in\Fp^3,\ t\in T\}.
\end{equation}
By the following lemma, $V'_T$ has a small intersection with every affine line and affine plane conic.

\begin{lemma}
\label{lem:avoidance}
Consider $V'_T\subseteq\Fp^4$ as defined in~\eqref{eq:selection}.
\begin{itemize}
\item[(a)]
Every affine line defined over $\Fp$ contains at most $3L$ points from
$V'_T$.
\item[(b)]
Every affine plane conic defined over $\Fp$ contains at most $6L$ points from
$V'_T$.
\end{itemize}
\end{lemma}

\begin{proof}
It suffices to consider $L=1$ since $V'_{T}=\bigcup_{t\in T}V'_{\{t\}}$.
Write $T=\{t\}$.

For (a), put $F_t(c,r):=N(c)+t-r$ so that $V'_{\{t\}}$ is the zero set of $F_t$. 
On an affine line parametrized by $(c_0,r_0)+s(v,w)$, the restriction of $F_t$ is
\[
    N(c_0+sv)+t-r_0-sw.
\]
Since $(v,w)\neq0$, we have that $v=0$ implies $w\ne0$, in which case this polynomial has degree~$1$.
Otherwise, it has degree~$3$ since its cubic coefficient is $N(v)\ne0$.
Thus every affine line meets $V'_{\{t\}}$ in at most three points.

For (b), let $C$ be a conic in an affine plane $\Pi$ defined over $\Fp$.
The restriction of $F_t$ to $\Pi$ has degree exactly~$3$ since $\Pi$ contains an affine line. 
If this plane cubic were reducible over $\Fp$, it would have a linear factor over $\Fp$. 
On the line defined by that factor, the restriction of $F_t$ would be the zero polynomial, contradicting our analysis in part~(a). 
Thus the cubic is irreducible over $\Fp$, and since $C$ is defined by a polynomial of degree~$2$, the cubic is relatively prime to this quadratic. 
Finally, the zero sets of these polynomials do not have a common component over $\overline{\Fp}$, since a common factor over $\overline{\Fp}$ would determine a common factor over $\Fp$ by multiplying its distinct Frobenius conjugates.
We may therefore apply B\'ezout's theorem, which gives $|C\cap V'_{\{t\}}|\le 3\cdot2=6$.
\end{proof}

\subsection{Putting everything together}

\begin{proposition}\label{prop:edges}
Every $U\subseteq V'_T$ of size $s\le p^2$ satisfies $e(U)\lesssim Ls^{3/2}$.
\end{proposition}

\begin{proof}
Partition $U=U_1\sqcup\cdots\sqcup U_L$ according to the translates in~\eqref{eq:selection}. 
Fix $i,j\in[L]$ and put $W=U_i\cup U_j$.
By Lemma~\ref{lem:avoidance}, each conic contains at most $12$~points
of $W$.
If $|W|>12^2$, then $W$ does not satisfy Theorem~\ref{thm:zahl}(b), and so we must have
$|\mathcal P(W)|\lesssim |W|^{3/2}$.
Next, Lemma~\ref{lem:avoidance} gives that each pencil contains at most $6$~points
of $W$, and therefore contributes at most $\binom{6}{2}$ edges. 
Thus $e(W)\leq|\mathcal{P}(W)|\cdot\binom{6}{2}\lesssim |W|^{3/2}$.
For $|W|\le 12^2$, the same bound follows from $e(W)\le\binom{|W|}{2}$.

To obtain the result, we apply the fact that every edge of $U$ resides in some $U_i\cup U_j$, along with the bounds $e(U_i\cup U_j)\lesssim |U_i\cup U_j|^{3/2}$ and $|U_i\cup U_j|\leq|U|=s$, and the fact that each vertex belongs to exactly $L$ of these unions:
\[
e(U)
\le\sum_{i\le j}e(U_i\cup U_j)
\lesssim 
\sum_{i\le j}|U_i\cup U_j|^{3/2}
\leq\sqrt{s}\sum_{i\le j}|U_i\cup U_j|
=Ls^{3/2}. \qedhere
\]
\end{proof}

\begin{proof}[Proof of Theorem~\ref{thm:main}]
Given a prime $p\equiv 3\bmod 4$, let $A$ denote the adjacency matrix of the strongly regular graph $\mathrm{VO}_4^-(p)$.
Taking $m$, $\mu$, and $\nu$ as in~\eqref{eq.params}, then $G:=I-\mu A+\nu\bar{A}$ is the Gram matrix of a two-distance tight frame $\Phi\in\mathbb{R}^{m\times p^4}$, which one can obtain from $G$ in polynomial time by Cholesky factorization.\footnote{This is standard, but for the sake of concreteness, we supply some details. First, we order the vertices lexicographically and compute the Cholesky factor $R$ of~$G$ so that $G=R^*R$.
Since $mG$ has integer entries of magnitude at most~$m$, its minors have polynomial bit length as a consequence of Hadamard's inequality, and so one may use rational symmetric elimination followed by square roots to compute $R$ in polynomial time to the requested precision.
After deleting the zero rows that belong to zero pivots, we obtain $\Phi$.}
Let $\Phi_T$ consist of the columns of $\Phi$ indexed by $V'_T$, defined in~\eqref{eq:selection}.\footnote{Note that the columns indexed by $V'_T$ can be listed in polynomial time by testing monic cubics for roots and computing $N(c)$ as the determinant of multiplication by~$c$.
Combined with the previous footnote, this implies that the matrix $\Phi_T$ is explicit.}
For $S\subseteq V'_T$ with $1\le |S|\le k\le p^2$, since the principal submatrix $G_S$ satisfies $|G_S-I|=\mu A_S+\nu\bar A_S$, the bound~\eqref{eq:Lev} and Gershgorin's circle theorem imply 
\[
\big\||G_S-I|\big\|_{2\to2}
\leq \mu\|A_S\|_{2\to2}+\nu\|\bar{A}_S\|_{2\to2}
\lesssim \mu \log(2k)\|A_S\|_{\DR}+\nu\|\bar{A}_S\|_{1\to1}.
\]
Next, Proposition~\ref{prop:edges} and Lemma~\ref{lem:density-norm} together give
$\|A_S\|_{\DR}\lesssim L\sqrt{k}$.
Since $\|\bar A_S\|_{1\to1}\le k$, we have
\begin{equation}\label{eq:final-RIP}
    \delta_k(\Phi_T)\le\Delta_k(\Phi_T)
    \lesssim\mu L\sqrt{k}\log(2k)+k\nu
    \sim \frac{L\sqrt{k}\log(2k)}{p}+\frac{k}{p^2}.
\end{equation}
Now let $p\to\infty$ through the primes congruent to $3$ modulo $4$, of which there are infinitely many, and take
\[
L=\lfloor\log p\rfloor,
\qquad 
T=\{0,\ldots,L-1\},
\qquad
n=|V'_T|=Lp^3,
\qquad
k=\left\lfloor\frac{p^2}{(\log p)^6}\right\rfloor.
\]
For all sufficiently large such $p$, we have $1\le L\le p$ and $1\le k\le p^2$. Then $m\to\infty$, $m/n\sim1/L\to0$, and $k=m^{2/3-o(1)}$, and by~\eqref{eq:final-RIP}, we also have $\delta_k(\Phi_T)\to0$.
\end{proof}

\section{The \texorpdfstring{$2/3$}{2/3} barrier and the role of phase cancelation}\label{sec:barrier}

In this paper, we constructed explicit restricted isometries with $k=m^{2/3-o(1)}$.
By Theorem~\ref{thm:absolute-barrier}, we have reached the limit of phase agnosticism. 
It remains to break the $2/3$~barrier:

\begin{problem}
\label{prob:beyond}
Construct an explicit family of $m\times n$ matrices satisfying
\[
m\longrightarrow\infty,
\qquad m/n\longrightarrow 0,
\qquad
\delta_k\longrightarrow 0,
\qquad 
k=m^{\alpha-o(1)}
\]
for some $\alpha>2/3$.
\end{problem}

It turns out that the matrices we used to prove Theorem~\ref{thm:main} do not break the $2/3$~barrier as an artifact of optimizing them for our proof technique:

\begin{proposition}\label{prop:actual-barrier}
Consider a sequence of matrices $\Phi\in\mathbb R^{m\times n}$, each a subframe of a real unit norm two-distance tight frame of coherence~$\mu$, such that $m\to\infty$ and $m/n\to0$.
Then for all sufficiently large $m$, it holds that
\[
\delta_k
\ge\min\big\{\mu\sqrt{k-1},1\big\}
\]
for every $1\le k\le n$.
In particular, if $\mu\gtrsim m^{-1/3}$,
then 
$\delta_k\to0$ requires $k=o(m^{2/3})$.
\end{proposition}

\begin{proof}
Let $\eta$ denote the smaller-magnitude inner product of the frame, so that $|\eta|=\nu\le\mu$. 
If the full frame has $N$ columns, then tightness gives a squared off-diagonal sum of $N/m-1$ in each column of its Gram matrix. 
Thus $(N-1)\nu^2\le N/m-1$, and so $\nu\le m^{-1/2}$.

Put $G=\Phi^*\Phi$, and let $D$ be the largest number of off-diagonal entries of magnitude $\mu$ in a column of~$G$. 
The matrix $H:=G-\eta J$ has rank at most $m+1$ and diagonal entries $1-\eta$. 
Each column of $H$ has at most $D$ nonzero off-diagonal entries, each of magnitude at most $2\mu$. 
Applying Cauchy--Schwarz to the nonzero eigenvalues of $H$ gives
\[
\frac{n(1-\eta)^2}{m+1}
=\frac{(\operatorname{tr}H)^2}{n(m+1)}
\le\frac{\|H\|_F^2}{n}
\le (1-\eta)^2+D(2\mu)^2.
\]
Since $\eta\to0$ and $n/m\to\infty$, it follows that $D\mu^2\to\infty$.

Choose a column index $i$ attaining $D$, and let $S\subseteq[n]$ consist of $i$ as well as the row indices of $\min\{k-1,D\}$ entries of magnitude $\mu$ in column $i$ of $G$. 
Then
\[
\delta_k^2\ge\|(G_S-I)e_i\|_2^2
=\min\{(k-1)\mu^2,D\mu^2\},
\]
which is eventually at least $\min\{(k-1)\mu^2,1\}$ since $D\mu^2\to\infty$.
\end{proof}

To solve Problem~\ref{prob:beyond}, one might instead apply the hybrid estimate from Section~\ref{sec:motivation}: retain the induced $1$-norm for the small entries, but estimate the \emph{phased} rectangular sums of the large entries using ideas from Bourgain et al.~\cite{BDFKK}.
We leave this for future work.

\section*{Acknowledgments}
CK was supported in part by NSF grant CIF 2549926. DGM was supported in part by NSF grant DMS 2220304 and AFOSR grant FA9550-26-1-B224.
DGM thanks OpenAI for complimentary access to ChatGPT-6 Pro.
More importantly, DGM thanks Matthew Fickus, Joey Iverson, John Jasper, and Emily King for over a decade of collaboration and friendship that helped shape his approach to digesting this solution.

\section*{AI and formalization declaration}

Theorem~\ref{thm:main} was proved in one shot using GPT-6~Astra.
While digesting the proof, the authors distilled the Gershgorin--flat RIP hybrid estimate, identified the 2/3~barrier, and explored the general approach of taking subframes of two-distance tight frames.
We conducted these investigations over the course of several long conversations with ChatGPT. To distinguish from AI slop, we took care to explain how the proof works, how one might arrive at it, and how its main ideas relate to prior work.
We used ChatGPT to generate an initial draft, and then, after we completely rewrote the paper, we used ChatGPT and Opus~5.5 to provide multiple rounds of feedback. 
All mathematical claims, citations, and computations were checked by the authors, who take full responsibility for the content of this paper.

All theorems, lemmas and propositions of this paper, including the case of Zahl’s theorem used here (Theorem~\ref{thm:zahl}),
have been formalized and machine-checked in Lean~4 with the Mathlib library; the running time of the construction in Theorem~\ref{thm:main} is not part of the formalization. The Lean formalization is maintained in the repository 
\begin{center}
\href{https://github.com/ckuemmerle/ExplicitRIPTwoThirds}{\textcolor{BrickRed}{\texttt{github.com/ckuemmerle/ExplicitRIPTwoThirds}}}
\end{center}
and was created by an Opus~5.5 orchestrator with GPT-6.1~Sol and GPT-6~Astra subagents.

\appendix

\section{Proof of Lemma~\ref{lem:srg-scaling}}

Taking a diagonal entry of $G^2=(\srg{v}/m)G$ and an entry of $G\one=0$ gives
\[
    1+\srg{k}\mu^2+(\srg{v}-1-\srg{k})\nu^2=\frac{\srg{v}}{m},
    \qquad
    \srg{k}(\mu+\nu)=(\srg{v}-1)\nu-\varepsilon.
\]
Collecting the terms involving $\srg{k}$ in the first identity and
substituting the second gives
\[
    \frac{\srg{v}}{m}-1
    = (\srg{v}-1)\nu^2+(\mu-\nu)\srg{k}(\mu+\nu)
    = (\srg{v}-1)\mu\nu-\varepsilon(\mu-\nu).
\]
Since $\srg{v}/m\to\infty$ and $\mu-\nu=o(1)$, we obtain $\mu\nu\sim1/m$,
and hence $\nu\asymp m^{-2/3}$. In particular, $\srg{v}\nu\to\infty$
and $\nu=o(\mu)$, so
\[
    \frac{\srg{k}}{\srg{v}}
    =\frac{(\srg{v}-1)\nu-\varepsilon}{\srg{v}(\mu+\nu)}
    \sim\frac{\nu}{\mu}\sim\frac{1}{m\mu^2}\asymp m^{-1/3}.
\]
Since $A$ acts by a scalar on the range of $G$, comparing diagonal
entries gives $AG=\varepsilon\srg{k}\mu G$. Comparing entries on a
nonedge and an edge, respectively, gives
\[
    (\mu+\nu)\srg{\mu}=\srg{k}\nu(1-\varepsilon\mu),\qquad
    (\mu+\nu)\srg{\lambda}
    =(\srg{k}-1)\nu+\varepsilon(\srg{k}\mu^2-1).
\]
The first identity gives
$\srg{\mu}\sim\srg{k}\nu/\mu\sim\srg{k}^2/\srg{v}$.
Dividing the second by the first gives
\[
    \frac{\srg{\lambda}}{\srg{\mu}}
    =\frac{1-1/\srg{k}+\varepsilon\bigl(\mu^2/\nu-1/(\srg{k}\nu)\bigr)}
           {1-\varepsilon\mu}
    =1+\varepsilon m\mu^3+o(1),
\]
since $\srg{k}\nu\asymp\srg{v}/m\to\infty$ and
$\mu^2/\nu\sim m\mu^3\asymp1$. This ratio tends to zero precisely when
$\varepsilon=-1$ eventually and $m\mu^3\to1$, in which case
$\mu\sim m^{-1/3}$ and $\nu\sim1/(m\mu)\sim m^{-2/3}$.

\end{document}